\documentclass[12pt,a4paper]{article}
\usepackage{indentfirst}
\usepackage{amsmath,amssymb,amsfonts,amsthm,graphics}
\usepackage{makeidx}
\usepackage{amsmath,amssymb,amsfonts,amsthm,graphics,graphicx}
\usepackage{makeidx}
\usepackage{color}
\usepackage{ifpdf}
\usepackage{subfigure}

\newtheorem{thm}{Theorem}

\newtheorem{pro}{Proposition}
\newtheorem{cor}{Corollary}

\theoremstyle{definition}

\def\-{\mbox{--}}

\def\pf{\noindent {\it Proof.} }
\begin{document}

\title{\Large\bf Total proper connection and graph operations}
\author{\small Yingying Zhang, Xiaoyu Zhu\\
\small Center for Combinatorics and LPMC\\
\small Nankai University, Tianjin 300071, China\\
\small E-mail: zyydlwyx@163.com; zhuxy@mail.nankai.edu.cn}
\date{}
\maketitle
\begin{abstract}

A graph is said to be {\it total-colored} if all the edges and vertices of the graph are colored. A path in a
total-colored graph is a {\it total proper path} if $(i)$ any two adjacent edges on the path differ in color,
$(ii)$ any two internal adjacent vertices on the path differ in color, and $(iii)$ any internal vertex of the path
differs in color from its incident edges on the path. A total-colored graph is called {\it total-proper connected}
if any two vertices of the graph are connected by a total proper path of the graph. For a connected graph $G$,
the {\it total proper connection number} of $G$, denoted by $tpc(G)$, is defined as the smallest number of colors
required to make $G$ total-proper connected. In this paper, we study the total proper connection number for the graph operations. We find that $3$ is the total proper connection number for the join, the lexicographic product and the strong product of nearly all graphs. Besides, we study three kinds of graphs with one factor to be traceable for the Cartesian product as well as the permutation graphs of the star and traceable graphs. The values of the total proper connection number
for these graphs are all $3$.


\end{abstract}

\section{Introduction}

In this paper, all graphs considered are simple, finite and undirected. We refer to the book \cite{B} for
undefined notation and terminology in graph theory. A path in an edge-colored graph is a {\it proper path}
if any two adjacent edges differ in color. An edge-colored graph is {\it proper connected} if any two
vertices of the graph are connected by a proper path of the graph. For a connected graph $G$, the
{\it proper connection number} of $G$, denoted by $pc(G)$, is defined as the smallest number of colors
required to make $G$ proper connected. Note that $pc(G)=1$ if and only if $G$ is a complete graph.
The concept of $pc(G)$ was first introduced by Borozan et al. \cite{BFG} and has been well-studied recently.
We refer the reader to \cite{ALL,BFG,GLQ,LLZ,MYWY1} for more details.

As a natural counterpart of the concept of proper connection, the concept of proper vertex connection
was introduced by Jiang et al. \cite{JLZZ}. A path in a vertex-colored graph is a {\it vertex-proper path}
if any two internal adjacent vertices on the path differ in color. A vertex-colored graph is
{\it proper vertex connected} if any two vertices of the graph are connected by a vertex-proper
path of the graph. For a connected graph $G$, the {\it proper vertex connection number} of $G$, denoted
by $pvc(G)$, is defined as the smallest number of colors required to make $G$ proper vertex connected.
Especially, set $pvc(G)=0$ for a complete graph $G$. Moreover, we have $pvc(G)\geq 1$ if $G$ is a
noncomplete graph.

Actually, the concepts of the proper connection and proper vertex connection were inspired by the concepts
of the rainbow connection and rainbow vertex connection. For details about them we refer to \cite{GMP,LiSS,LiS,MYWY}. Here we only state the concept of the total rainbow connection of graphs, which was
introduced by Liu et al. \cite{LMS} and also studied in \cite{JLZ, S}. A graph is {\it total-colored} if all the
edges and vertices of the graph are colored. A path in a total-colored graph is a {\it total
rainbow path} if all the edges and internal vertices on the path differ in color. A total-colored graph is
{\it total-rainbow connected} if any two vertices of the graph are connected by a total rainbow path of the graph.
For a connected graph $G$, the {\it total rainbow connection number} of $G$, denoted by $trc(G)$, is defined as the
smallest number of colors required to make $G$ total-rainbow connected. Motivated by the concept of the total rainbow
connection, for the proper connection and proper vertex connection Jiang et al. \cite{JLZ1} introduced the concept of the total proper
connection. A path in a total-colored graph is a {\it total proper path} if $(i)$ any two adjacent edges on the path
differ in color, $(ii)$ any two internal adjacent vertices on the path differ in color, and $(iii)$ any internal vertex of
the path differs in color from its incident edges on the path. A total-colored graph is {\it total-proper connected}
if any two vertices of the graph are connected by a total proper path of the graph. For a connected graph $G$,
the {\it total proper connection number} of $G$, denoted by $tpc(G)$, is defined as the smallest number of colors
required to make $G$ total-proper connected. It is easy to obtain that $tpc(G)=1$ if and only if $G$ is a complete
graph, and $tpc(G)\geq 3$ if $G$ is not complete. Moreover,
$$tpc(G)\geq \max\{pc(G),pvc(G)\}.\ \ \ \  \ \ \ \ \ \ \ \ \ \  \ \ \ \ \ \ \ \ (*)$$

We recall some fundamental results on $tpc(G)$ which can be found in \cite{JLZ1}.

\begin{pro}\label{pro1}\cite{JLZ1} If $G$ is a nontrivial connected graph and $H$ is a connected spanning subgraph of $G$, then $tpc(G)\leq tpc(H)$.
In particular, $tpc(G)\leq tpc(T)$ for every spanning tree $T$ of $G$.
\end{pro}

\begin{pro}\label{pro2}\cite{JLZ1} Let $G$ be a connected graph of order $n\geq 3$ that contains a bridge. If $b$ is the maximum number of bridges
incident with a single vertex in $G$, then $tpc(G)\geq b+1$.
\end{pro}

Let $\Delta(G)$ denote the maximum degree of a connected graph $G$. We have the following.

\begin{thm}\label{thm1}\cite{JLZ1} If $T$ is a tree of order $n\geq 3$, then $tpc(T)=\Delta(T)+1$.
\end{thm}

The consequence below is immediate from Proposition \ref{pro1} and Theorem \ref{thm1}.

\begin{cor}\label{cor1}\cite{JLZ1} For a nontrivial connected graph $G$,
$$tpc(G)\leq \min\{\Delta(T)+1:\ T\ is\ a\ spanning\ tree\ of\ G\}.$$
\end{cor}

A {\it Hamiltonian path} in a graph $G$ is a path containing every
vertex of $G$ and a graph having a Hamiltonian path is a {\it traceable graph}. We get the following result.

\begin{cor}\label{cor2}\cite{JLZ1} If $G$ is a traceable graph that is not complete, then $tpc(G)=3$.
\end{cor}

Let $K_{m,n}$ denote a complete bipartite graph, where $1\leq m\leq n$. Clearly,
$tpc(K_{1,1})=1$ and $tpc(K_{1,n})=n+1$ if $n\geq 2$. For $m\geq 2$, we have the result below.

\begin{thm}\label{thm12}\cite{JLZ1} For $2\leq m\leq n$, we have $tpc(K_{m,n})=3$.
\end{thm}

\begin{thm}\label{thm2}\cite{JLZ1} Let $G$ be a $2$-connected graph. Then $tpc(G)\leq 4$ and the upper bound is sharp.
\end{thm}

The standard products (Cartesian, direct, strong, and lexicographic) draw a constant attention of graph research community, see some papers \cite{ACKP,BS,GV,KS,NS,P,S1,Z}. In this paper we consider the join, permutation graph and three standard products: the Cartesian, the strong, and the lexicographic with respect to the total proper connection number. Each of them will be treated in one of the forthcoming sections. In Section 2,5 and 6, we prove that $3$ is the total proper connection number for the join, the lexicographic product and the strong product of nearly all graphs, respectively. In Section 3, we study three kinds of graphs with one factor to be traceable for the Cartesian product, and obtain that the values of the total proper connection number of these graphs are all $3$. In Section 4, we show that the total proper connection numbers of the permutation graphs of the star and traceable graphs are also $3$.

\section{Joins of graphs}

The {\it join} $G\vee H$ of two graphs $G$ and $H$ has vertex set $V(G)\cup V(H)$ and its edge set consists of $E(G)\cup E(H)$ and the set $\{uv:u\in V(G)$ and $v\in V(H)\}$.

\begin{thm}\label{thmJ1} If $G$ and $H$ are connected graphs such that $G\vee H$ is not complete, then $tpc(G\vee H)=3$.
\end{thm}

\pf If $G$ and $H$ are both nontrivial connected graphs such that $G\vee H$ is not complete, then $G\vee H$ contains the graph in Theorem \ref{thm12} as a spanning subgraph. By Proposition \ref{pro1} and Theorem \ref{thm12}, it follows that $tpc(G\vee H)=3$. Thus we may assume that $G$ is a nontrivial connected graph of order at least 3 that is not complete and $H=K_1$ where $V(H)=\{w\}$. Since $G\vee K_1$ is not complete, it follows that $tpc(G\vee K_1)\geq3$ and so it remains to show that $tpc(G\vee K_1)\leq3$. Let $T$ be a spanning tree of $G$. By Proposition \ref{pro1}, it suffices to show that $tpc(T\vee K_1)\leq3$. For a vertex $v\in V(T)$, let
$e_{T}(v)$ denote the eccentricity of $v$ in $T$, i.e., the maximum
of the distances between $v$ and the other vertices in $T$. Let
$V_{i}=\{u\in V(T) : d_T(u,v)=i\}$, where $0\leq i\leq e_{T}(v)$.
Hence $V_{0}=\{v\}$. Define a 3-coloring $c$ of the vertices and edges of $T\vee K_1$ by

\begin{eqnarray}c(x)=
\begin{cases}
1 &if\ x\in V_i,\ i\ is\ even\ and\ 0\leq i\leq e_{T}(v) \cr
2 &if\ x\in V_i,\ i\ is\ odd\ and\ 1\leq i\leq e_{T}(v) \cr
3 &if\ x=w;
\end{cases}
\end{eqnarray}

\begin{eqnarray}c(wx)=
\begin{cases}
1 &if\ x\in V_i,\ i\ is\ odd\ and\ 1\leq i\leq e_{T}(v) \cr
2 &if\ x\in V_i,\ i\ is\ even\ and\ 0\leq i\leq e_{T}(v);
\end{cases}
\end{eqnarray}

\begin{eqnarray}c(xy)=
\begin{cases}
3 &if\ x\in V_i,\ y\in V_{i+1}\ and\ 0\leq i\leq e_{T}(v)-1.
\end{cases}
\end{eqnarray}

Let $x$ and $y$ be two vertices of $T\vee K_1$. Since $w$ is adjacent to every vertex in $T$, we may assume $x,y\in V(T)$. First, suppose that $x\in V_i$ and $y\in V_j$, where $0\leq i<j\leq e_{T}(v)$. If $i$ and $j$ are of opposite parity, then the path $xwy$ is a total proper $x$-$y$ path in $T\vee K_1$. Thus we may assume that $i$ and $j$ are of the same parity and so $j-i\geq2$. Let $z\in V_{j-1}$ such that $yz$ is an edge of $T$. Then the path $xwzy$ is a total proper $x$-$y$ path in $T\vee K_1$. Next, suppose that $x,y\in V_i$ for some $i$ with $1\leq i\leq e_{T}(v)$. Let $z\in V_{i-1}$ such that $xz$ is an edge of $T$. Then the path $xzwy$ is a total proper $x$-$y$ path in $T\vee K_1$. Hence for any two vertices $x$ and $y$ in $T\vee K_1$, there exists a total proper path between them and so $tpc(T\vee K_1)\leq3$. Therefore, $tpc(G\vee H)=3$.\qed

\section{The Cartesian product}

The {\it Cartesian product} $G\square H$ of two graphs $G$ and $H$ is the graph with vertex set $V(G)\times V(H)$, in which two vertices $(g,h)$ and $(g',h')$ are adjacent if and only if $g=g'$ and $hh'\in E(H)$, or $h=h'$ and $gg'\in E(G)$. Clearly, the Cartesian product is commutative, that is, $G\square H$ is isomorphic to $H\square G$. Moreover, $G\square H$ is $2$-connected whenever $G$ and $H$ are connected. From Theorem \ref{thm2}, we have that $3\leq tpc(G\square H)\leq4$. In this section, we mainly study three kinds of graphs with one factor to be traceable for the Cartesian product, and obtain that the values of the total proper connection number of these graphs are all $3$.

\begin{thm}\label{thmC1} Let $G$ and $H$ be two nontrivial traceable graphs with $|G|=n$ and $|H|=m$. Then $tpc(G\square H)=3$.
\end{thm}

\pf Clearly, $P_n$ and $P_m$ are spanning subgraphs of $G$ and $H$, respectively. Since $P_n\square P_m$ is traceable, it follows that $tpc(P_n\square P_m)=3$ by Corollary \ref{cor2}. Moreover, $P_n\square P_m$ is a spanning subgraph of $G\square H$. From Proposition \ref{pro1}, we have that $tpc(G\square H)\leq tpc(P_n\square P_m)$. Thus $tpc(G\square H)\leq3$. Since $G\square H$ is not complete, $tpc(G\square H)\geq3$ and so $tpc(G\square H)=3$.\qed

\begin{thm}\label{thmC2} Let $G$ be a nontrivial traceable graph and $H$ be a connected graph with maximum degree $|H|-1$. Then $tpc(G\square H)=3$.
\end{thm}

\pf Since $G\square H$ is not complete, we just need to show that $tpc(G\square H)\leq3$. Let $P_n=g_1g_2...g_n$ be a spanning subgraph of the nontrivial traceable graph $G$, where $n\geq2$. And let $K_{1,s}$ be a spanning subgraph of $H$ with $V(K_{1,s})=\{h_0,h_1,...,h_s\}$, where $s=|H|-1$ and $h_0$ is the central vertex. Then $P_n\square K_{1,s}$ is a spanning subgraph of $G\square H$ and so it suffices to show that $tpc(P_n\square K_{1,s})\leq3$ by Proposition \ref{pro1}. From Theorem \ref{thmC1}, we only need to consider the case that $s\geq3$.

Define a 3-coloring $c$ of the vertices and edges of $P_n\square K_{1,s}$ by

\begin{eqnarray}c(g_i,h_j)=
\begin{cases}
1 &if\ i\in [n]\ and\ j=0, \cr
2 &if\ either\ i\in [n],\ i\ is\ odd\ and\ 2\leq j\leq s\\
&\ \ \ or\ i\in [n],\ i\ is\ even\ and\ j=1, \cr
3 &otherwise;
\end{cases}
\end{eqnarray}

\begin{eqnarray}c((g_i,h_j)(g_{i+1},h_j))=
\begin{cases}
1 &if\ 1\leq i\leq n-1\ and\ j\in[s], \cr
2 &if\ 1\leq i\leq n-1\ and\ j=0;
\end{cases}
\end{eqnarray}

\begin{eqnarray}c((g_i,h_0)(g_i,h_j))=
\begin{cases}
2 &if\ either\ i\in [n],\ i\ is\ odd\ and\ j=1\\
&\ \ \ or\ i\in [n],\ i\ is\ even\ and\ 2\leq j\leq s, \cr
3 &otherwise.
\end{cases}
\end{eqnarray}

It remains to check that there is a total proper path between any two vertices $(g_i,h_k),(g_j,h_t)$ in $P_n\square K_{1,s}$, where $1\leq i,j\leq n$ and $0\leq k,t\leq s$. For $i=j$, if $k=0$ or $t=0$, then the edge $(g_i,h_k)(g_j,h_t)$ is the desired path; if $k=1$ or $t=1$, then the desired path is $(g_i,h_k)(g_i,h_0)(g_j,h_t)$; if $2\leq k,t\leq s$, then the desired path is $(g_i,h_k)(g_i,h_0)(g_i,h_1)(g^*,h_1)\\(g^*,h_0)(g^*,h_t)(g_j,h_t)$, where $g^*$ is a neighbor of $g_i$ in $P_n$. For $2\leq i_1,i_p\leq s$, where $2\leq p\leq n$ and $p$ is even, set $P=(g_1,h_{i_1})(g_1,h_0)(g_1,h_1)(g_2,h_1)(g_2,h_0)(g_2,h_{i_2})(g_3,h_{i_2})(g_3,h_0)(g_3,h_1)\\
(g_4,h_1)(g_4,h_0)(g_4,h_{i_4})
...(g_r,h_{i_{r-1}})(g_r,h_0)(g_r,h_1)(g_{r+1},h_1)(g_{r+1},h_0)(g_{r+1},h_{i_{r+1}})
...(g_{n-1},h_{i_{n-2}})\\(g_{n-1},h_0)(g_{n-1},h_1)(g_n,h_1)(g_n,h_0)(g_n,h_{i_n})$ when $n$ is even and $P=(g_1,h_{i_1})(g_1,h_0)(g_1,h_1)\\(g_2,h_1)(g_2,h_0)(g_2,h_{i_2})(g_3,h_{i_2})(g_3,h_0)(g_3,h_1)(g_4,h_1)(g_4,h_0)(g_4,h_{i_4})
...(g_r,h_{i_{r-1}})(g_r,h_0)(g_r,h_1)\\(g_{r+1},h_1)(g_{r+1},h_0)(g_{r+1},h_{i_{r+1}})
...(g_{n-2},h_{i_{n-3}})(g_{n-2},h_0)(g_{n-2},h_1)(g_{n-1},h_1)(g_{n-1},h_0)(g_{n-1},h_{i_{n-1}})\\(g_n,h_{i_{n-1}})(g_n,h_0)(g_n,h_1)$ when $n$ is odd. According to the total-coloring $c$ of $P_n\square K_{1,s}$, the path $P$ is total proper. For $i\neq j$, we can always find a total proper path which is a subpath of $P$ between $(g_i,h_k)$ and $(g_j,h_t)$ in $P_n\square K_{1,s}$. Thus $tpc(P_n\square K_{1,s})\leq3$ and the proof is complete.\qed

\begin{thm}\label{thmC3} Let $G$ be a nontrivial traceable graph and $H$ be a connected graph with maximum degree $|H|-2$. Then $tpc(G\square H)=3$.
\end{thm}

\pf If $|H|=4$, then $H$ is traceable and so $tpc(G\square H)=3$ by Theorem \ref{thmC1}. Thus we only need to consider the case that $|H|\geq5$. Since $diam(G\square H)\geq2$, we have $tpc(G\square H)\geq3$ and so it remains to show that $tpc(G\square H)\leq3$. Let $P_n=g_1g_2...g_n$ be a spanning subgraph of the nontrivial traceable graph $G$. Denote by $x$ a vertex of $H$ with the maximum degree $|H|-2$ and $z$ the unique vertex not adjacent to $x$ in $H$. Since $H$ is connected, $z$ must be adjacent to one neighbor, say $y$, of $x$. We then take a spanning tree $T$ of $H$ containing the edges $zy$ and $xu$, where $u\in V(H)\backslash\{x,z\}$. Clearly, $P_n\square T$ is a spanning subgraph of $G\square H$. From Proposition \ref{pro1}, it suffices to show that $tpc(P_n\square T)\leq3$.

First suppose that $n=2$. Define a $3$-coloring $c$ of the vertices and edges of $P_2\square T$ as follows. For $w\in V(H)\backslash\{x,y,z\}$, set $c(g_1,y)=c((g_1,x)(g_1,w))=c(g_2,w)=c((g_2,x)(g_2,y))=1$, $c((g_1,x)(g_1,y))=c(g_1,w)=c((g_2,x)(g_2,w))=c(g_2,y)=2$ and $c(g_1,x)=c((g_1,w)(g_2,w))=c(g_2,x)=c((g_1,y)(g_2,y))=c((g_1,z)(g_1,y))=c((g_2,z)(g_2,y))=3$. Moreover, give each of the unmentioned vertices and edges in $P_2\square T$ a random color from $\{1,2,3\}$. Next it remains to check that there is a total proper path between any two vertices $(g_i,h),(g_j,h')$ in $P_2\square T$. According to the total-coloring $c$ of $P_2\square T$, it is easy to see that the path $P=(g_1,z)(g_1,y)(g_1,x)(g_1,w)(g_2,w)(g_2,x)(g_2,y)(g_2,z)$ is total proper. For $i=j$, if $h,h'\in V(T)\backslash\{x,y,z\}$, then the path $(g_i,h)(g_i,x)(g_i,y)(g_{3-i},y)(g_{3-i},x)(g_{3-i},h')(g_j,h')$ is the desired path; otherwise, we can find a total proper path which is a subpath of $P$ between $(g_i,h)$ and $(g_j,h')$. For $i\neq j$, if $h=h'$, then the edge $(g_i,h)(g_j,h')$ is the desired path; if $h,h'\in V(T)\backslash\{x,y,z\}$, the total proper path is $(g_i,h)(g_i,x)(g_i,y)(g_j,y)(g_j,x)(g_j,h')$; otherwise, we can always find a total proper path which is a subpath of $P$ between $(g_i,h)$ and $(g_j,h')$. Thus $tpc(P_2\square T)\leq3$.

Then suppose that $n=3$. On the basis of the total-coloring $c$ we give above, color the vertices and edges of $P_3\square T$ in such a way that for any $w\in V(G)\backslash\{x,y,z\}$, the trail $(g_3,w)(g_3,x)(g_3,y)(g_3,z)(g_2,z)(g_1,z)(g_1,y)(g_1,x)(g_1,w)(g_2,w)(g_3,w)$ is total-proper connected. Again for the remaining edges and vertices of $P_3\square T$, give them any color from $\{1,2,3\}$ as you like. Similar to the above checking process, we can get that there is a total proper path between any two vertices $(g_i,h),(g_j,h')$ in $P_3\square T$ and so $tpc(P_3\square T)\leq3$.

Finally suppose that $n\geq4$. We divide our discussion into three cases:

\textbf{Case 1.} $n\equiv1\ (\text{mod}\ 3)$.

We give a total-coloring of $P_n\square T$ using the color set $\{1,2,3\}$ in such a way that for any $w\in V(T)\backslash\{x,y,z\}$ and $2\leq i\leq n$, the trail $(g_i,w)(g_i,x)(g_i,y)(g_i,z)(g_{i-1},z)\cdots (g_1,z)(g_1,y)\\(g_1,x)(g_1,w)(g_2,w)\cdots (g_n,w)$ is total-proper connected. Since $n\equiv1\ (\text{mod}\ 3)$, we have $l(P_n)\equiv0\ (\text{mod}\ 3)$. Moreover, $diam(T)=3$. Thus it is easy to find that the path $(g_1,w)(g_2,w)\cdots(g_n,w)(g_n,x)(g_n,y)(g_n,z)$ is also total proper. So far we have confirmed the colors of all the vertices and some edges of $P_n\square T$. For the other uncolored edges, which of course, are all in form of $(g_i,h)(g_j,h)$, give this kind of edge a color differing from the colors which its endpoints have already used. Thus we can check that for $2\leq i\leq n-1$, the paths $(g_n,w)(g_n,x)(g_{n-1},x)\cdots(g_i,x)$ and $(g_n,w)(g_n,x)(g_n,y)(g_{n-1},y)\cdots(g_i,y)$ are all total proper. Next it remains to show that there is a total proper path between any two vertices $(g_i,h),(g_j,h')$ in $P_n\square T$. For $i=j$, if $h,h'\in V(T)\backslash\{x,y,z\}$, then the desired path is $(g_1,h)(g_2,h)...(g_n,h)(g_n,x)(g_n,y)(g_n,z)(g_{n-1},z)...(g_1,z)(g_1,y)(g_1,x)(g_1,h')$ when $i=j=1$ and the desired path is $(g_i,h)(g_i,x)(g_i,y)(g_i,z)(g_{i-1},z)...(g_1,z)(g_1,y)(g_1,x)(g_1,h')(g_2,h')\\...(g_j,h')$ when $i=j\geq2$; otherwise, we can find a total proper $(g_i,h)$-$(g_j,h')$ path which is a subpath of $(g_i,w)(g_i,x)(g_i,y)(g_i,z)$, where $w\in V(T)\backslash\{x,y,z\}$. Now we assume that $i\neq j$, say $i< j$. For $i=1$, the path $(g_1,h)...(g_1,z)(g_2,z)...(g_j,z)...(g_j,h')$ is the desired path. For $i\geq2$, if $h=h'$, then the path $(g_i,h)(g_{i+1},h)...(g_j,h')$ is the desired path; otherwise, the desired path is $(g_i,h)...(g_i,z)(g_{i-1},z)...(g_1,z)(g_1,y)(g_1,x)(g_1,h')(g_2,h')...(g_j,h')$ when $h'\in V(T)\backslash\{x,y,z\}$ and $(g_i,h)...(g_i,z)(g_{i-1},z)...(g_1,z)(g_1,y)(g_1,x)(g_1,w)(g_2,w)...(g_n,w)...\\(g_n,h')(g_{n-1},h')\cdots (g_j,h')$ when $h'\in\{x,y,z\}$, where $w\in V(T)\backslash\{x,y,z,h\}$. Thus $tpc(P_n\square T)\\\leq3$.

\textbf{Case 2.} $n\equiv2\ (\text{mod}\ 3)$.

This case can be viewed as adding one $T$-layer to the graph in Case 1. So we give the coloring in Case 1 to $P_n\square T$ except for the last $T$-layer, which is the $Z$-induced subgraph where $Z=\{(g_n,v):v\in V(T)\}$. We color this induced subgraph in such a way that for any $w\in V(T)\backslash\{x,y,z\}$, the trail $(g_n,w)(g_n,x)(g_n,y)(g_n,z)(g_{n-1},z)\cdots(g_1,z)(g_1,y)(g_1,x)(g_1,w)\\(g_2,w)\cdots(g_n,w)$ is total-proper connected. Similar to the checking process in Case 1, we can check that for any two vertices $(g_i,h),(g_j,h')$ in $P_n\square T$, there is a total proper path between them and so $tpc(P_n\square T)\leq3$.

\textbf{Case 3.} $n\equiv0\ (\text{mod}\ 3)$.

The last case again can be viewed as one $T$-layer more added to Case 2. We likewise color the former $n-1$ $T$-layers as we have discussed in Case 2 and for the last $T$-layer, again make the trail $(g_n,w)(g_n,x)(g_n,y)(g_n,z)(g_{n-1},z)\cdots(g_1,z)(g_1,y)(g_1,x)(g_1,w)(g_2,w)\cdots(g_n,w)$ total-proper connected, where $w\in V(T)\backslash\{x,y,z\}$. Moreover, color the two edges $(g_{n-1},x)(g_{n-2},x)$ and $(g_{n-1},y)(g_{n-2},y)$ with the colors of the edges $(g_{n-2},x)(g_{n-2},y)$ and $(g_{n-2},y)(g_{n-2},z)$, respectively. Next it remains to show that there is a total proper path between any two vertices $(g_i,h),(g_j,h')$ in $P_n\square T$. By symmetry, suppose that $i=n$ and $j=n-1$. If $h'=z$, then the path $(g_n,h)...(g_n,z)(g_{n-1},z)$ is the desired path; otherwise, the desired path is $(g_n,h)...(g_n,z)(g_{n-1},z)...(g_1,z)(g_1,y)(g_1,x)(g_1,w)(g_2,w)...(g_{n-2},w)...(g_{n-2},h')(g_{n-1},h')$, \\where $w\in V(T)\backslash\{x,y,z\}$. For the other cases, we can check in a similar way as Case 2. Thus $tpc(P_n\square T)\leq3$.\qed

\section{Permutation graphs}

Let $G$ be a graph with $V(G)=\{v_1,...,v_n\}$ and $\alpha$ be a permutation of $[n]$. Let $G'$ be a copy of $G$ with vertices labeled $\{u_1,...,u_n\}$ where $u_i\in G'$ corresponds to $v_i\in G$. Then the {\it permutation graph} $P_{\alpha}(G)$ of $G$ can be obtained from $G\cup G'$ by adding all edges of the form $v_iu_{\alpha(i)}$. This concept was first introduced by Chartrand and Harary \cite{CH}. Note that if $\alpha$ is the identity permutation on $[n]$, then $P_{\alpha}(G)=G\square K_2$ is the Cartesian product of a graph $G$ and $K_2$. Moreover, $P_{\alpha}(G)$ is $2$-connected whenever $G$ is connected. From Theorem \ref{thm2}, we have that $3\leq P_{\alpha}(G)\leq4$. In this section, we mainly study the permutation graphs of the star and traceable graphs, and obtain that the values of the total proper connection number of these graphs are all $3$.

\begin{thm}\label{thmP1} Let $G$ be a nontrivial traceable graph of order $n$. Then $tpc(P_{\alpha}(G))=3$ for each permutation $\alpha$ of $[n]$.
\end{thm}

\pf Let $P=v_1v_2...v_n$ be a hamiltonian path of $G$. Then $P'=u_1u_2...u_n$ is a hamiltonian path of $G'$. Besides, we write $P^{-1}$ and $P'^{-1}$ as the reverse of $P$ and $P'$, respectively. If $\alpha(n)=1$ or $n$, then clearly $P_{\alpha}(G)$ is traceable and the theorem holds according to Corollary \ref{cor2}. Otherwise, we suppose that $\alpha(n)=i\ (2\leq i\leq n-1)$. Since $P_{\alpha}(G)$ is not complete, it remains to show that $tpc(P_{\alpha}(G))\leq3$. Define a $3$-coloring $c$ of the vertices and edges of $P_{\alpha}(G)$ as follows. First color the vertices and edges of the path $P$ starting from $v_1$ in turn with the colors $1,2,3$. Then color the remaining vertices and edges in the three paths $v_1Pv_nu_iP'^{-1}u_1,v_1Pv_nu_iP'u_n$ and $u_{\alpha(1)}v_1Pv_n$ so that each follows the sequence $1,2,3,...,1,2,3,...$. Finally set $c(v_ju_{\alpha(j)})=c(v_{j-1}v_j)$, where $2\leq j\leq n-1$. Next we check that there is a total proper path between any two vertices in $P_{\alpha}(G)$. It is easy to see the total proper paths between all pairs of vertices except between $u_s$ and $u_t$ with $1\leq s\leq i-1$ and $i+1\leq t\leq n$. In this case, the path $u_sP'u_iv_nP^{-1}v_{\alpha^{-1}(t)}u_t$ is the desired total proper path. Thus the proof is complete.\qed

\begin{thm}\label{thmP2} Every permutation graph of a star of order at least $4$ has total proper connection number $3$.
\end{thm}

\pf For an integer $m\geq3$, let $G=K_{1,m}$ be the star with vertex set $\{v_0,v_1,...,v_m\}$, where $v_0$ is the central vertex. Then there are exactly two non-isomorphic permutation graphs, namely $P_{\alpha_1}(G)=G\square K_2$ where $\alpha_1$ is the identity permutation on the set $\{0,1,...,m\}$ and $P_{\alpha_2}(G)$ where $\alpha_2=(0,1)$. By Theorem \ref{thmC2}, we have that $tpc(P_{\alpha_1}(G))=3$. It remains to show that $tpc(P_{\alpha_2}(G))=3$. Let $\{v'_0,v'_1,...,v'_m\}$ be the corresponding vertex set in the second copy $G'$ of $G$. Since $P_{\alpha_2}(G)$ is not complete, we just need to show that $tpc(P_{\alpha_2}(G))\leq3$.

Define a total-coloring $c$ of $P_{\alpha_2}(G)$ with three colors by assigning $(1)$ the color 1 to the vertices $v_0,v'_0$ and the edges $v_iv'_i$ for $2\leq i\leq m$, $(2)$ the color 2 to the vertices $v'_2,v_i$ and the edges $v_0v'_1,v_0v_2,v'_0v'_i$ for $i\in[m]\backslash\{2\}$ and $(3)$ the color 3 to the remaining vertices and edges of $P_{\alpha_2}(G)$. It remains to check that there is a total proper path between any two vertices $u,v$ in $P_{\alpha_2}(G)$. For $3\leq i\leq m$, the cycle $C_i=v_0v_2v'_2v'_0v'_iv_iv_0$ is a total-proper connected $6$-cycle. Thus we may assume that $u$ and $v$ do not belong to any one of the $m-2$ cycles at the same time.

First suppose that $u=v_1$ and $v=v'_1$ by symmetry. Then the path $uv_0v$ or $uv'_0v$ is the desired path. Next suppose that $u=v_1$ or $v'_1$ and $v\in V(P_{\alpha_2}(G))\backslash\{v_1,v'_1\}$ by symmetry. If $v=v_0$ or $v'_0$, then the edge $uv$ is the desired path. Now assume first that $u=v_1$. If $v=v_2$, then $uv_0v$ is the desired path, while if $v=v'_2$, then $uv_0v_2v$ is the desired path. For $i\geq3$, if $v=v_i$, then $uv'_0v'_iv$ is the desired path, while if $v=v'_i$, then $uv'_0v$ is the desired path. Then assume that $u=v'_1$. If $v=v_2$, then $uv'_0v'_2v$ is the desired path, while if $v=v'_2$, then $uv'_0v$ is the desired path. For $i\geq3$, if $v=v_i$, then $uv_0v$ is the desired path, while if $v=v'_i$, then $uv_0v_iv$ is the desired path. Finally suppose that $u,v\in V(P_{\alpha_2}(G))\backslash\{v_0,v'_0,v_1,v'_1,v_2,v'_2\}$. Let $P=v_iv'_iv'_0v'_2v_2v_0v_jv'_j$, where $3\leq i,j\leq m$ and $i\neq j$. According to the total coloring $c$, it is easy to see that the path $P$ is a total proper path. Moreover, we can always find a total proper path which is a subpath of $P$ between $u$ and $v$. Thus $tpc(P_{\alpha_2}(G))\leq3$ and we complete the proof.\qed

We conclude this section with the following question: Is there a class of nontrivial connected graphs $G$ such that $tpc(P_{\alpha}(G))=4$ for some permutation graph $P_{\alpha}(G)$ of $G$?

\section{The lexicographic product}

The {\it lexicographic product} $G\circ H$ of graphs $G$ and $H$ is the graph with vertex set $V(G)\times V(H)$, in which two vertices $(g,h),(g',h')$ are adjacent if and only if $gg'\in E(G)$, or $g=g'$ and $hh'\in E(H)$. The lexicographic product is not commutative and is connected whenever $G$ is connected. In a tree $T$, we denote the parent of the vertex $v$ by $p(v)$.

\begin{thm}\label{lexicographic product} Let $G$ and $H$ be two nontrivial graphs. If $G$ is connected and $G\circ H$ is not complete, then $tpc(G\circ H)=3$.
\end{thm}

\pf Since $G\circ H$ is not complete, it follows that $tpc(G\circ H)\geq3$ and so we just need to show that $tpc(G\circ H)\leq3$. If $G$ has only two vertices, i.e. $G=K_2$, then $G\circ H$ contains the graph in Theorem \ref{thm12} as a spanning subgraph and so $tpc(G\circ H)\leq3$ by Proposition \ref{pro1} and Theorem \ref{thm12}. Now we may assume that $G$ is a nontrivial connected graph of order at least 3. Take a spanning tree $T$ from $G$ and appoint a pendant vertex of $T$, say $r$, to be the root of $T$. Since $r$ is a pendant vertex, it has only one neighbor in $T$ called $t$. For the graph $H$, we view it as an empty graph. Thus the lexicographic product $T\circ H$ is a spanning subgraph of $G\circ H$. By Proposition \ref{pro1}, it suffices to show that $tpc(T\circ H)\leq3$.

Define a total-coloring $c$ of $T\circ H$ using the color set $A=\{1,2,3\}$ as follows. Let $V(H)=\{h_1,h_2,\cdots,h_n\}$ and then set $X=\{(g,h_1)|g\in V(T)\}$. We first give the vertices and edges of $X$-induced subgraph of $T\circ H$ a total-coloring using $A$ in such a way that for any vertex $g\in V(T)$, the path $(g,h_1)(g_1,h_1)(g_2,h_1)\cdots(t,h_1)(r,h_1)$ in $T\circ H$ is total proper, where $gg_1g_2\cdots tr$ is the unique path between $g$ and $r$ in $T$. Then color the edge $(r,h_1)(t,h_2)$ in such a way that the path $(t,h_1)(r,h_1)(t,h_2)$ is total proper. Let $Y=\{(q,h_2)|q\in V(T)\backslash\{r\}\}$. We give the $Y$-induced subgraph of $T\circ H$ a total-coloring in such a way that for any two vertices $(g,h_1),(g',h_2)$ in $T\circ H$, the path $(g,h_1)(g_1,h_1)(g_2,h_1)\cdots(t,h_1)(r,h_1)(t,h_2)\cdots(g'_2,h_2)(g'_1,h_2)(g',h_2)$ is total proper, where $gg_1g_2\cdots tr$ and $g'g'_1g'_2\cdots tr$ are the paths from $g$ to $r$ and $g'$ to $r$ in $T$ respectively. For $(g,h_i)\in V(T\circ H)$, where $g\in V(T)\backslash\{r,t\}$ and $i\in[n]$, set $c((g,h_i)(p(g),h_1))=c((g,h_1)(p(g),h_1))$ and $c((g,h_i)(p(g),h_2))=c((g,h_2)(p(g),h_2))$. By the way, we let $c((t,h_1)(r,h_i))=c((t,h_1)(r,h_1))$ for $2\leq i\leq n$, $c((r,h_1)(t,h_j))=c((t,h_2)(r,h_j))=c((r,h_1)(t,h_2))$ for $3\leq j\leq n$, $c((r,h_3)(t,h_s))=c(t,h_2)$ for $4\leq s\leq n$ and $c((r,h_2)(t,h_2))=c(r,h_1)$. Pick one neighbor of the vertex $t$ in $T$ other than $r$ called $a$ and make $c((a,h_2)(t,h_i))\\=c(t,h_2)$ for $3\leq i\leq n$, $c(r,h_3)=c(t,h_3)=c((a,h_2)(t,h_2))$, and $c((t,h_3)(r,h_j))=c(a,h_2)$ for $3\leq j\leq n$.

Next it remains to check that for any two vertices $(g,h_i),(g',h_j)$ in $T\circ H$, where $g,g'\in V(T)$ and $h_i,h_j\in V(H)$, there is a total proper path between them. Let $gg_1g_2\cdots tr$ and $g'g'_1g'_2\cdots tr$ be the paths from $g$ to $r$ and $g'$ to $r$ in $T$, respectively. If $g,g'\in V(T)\backslash\{r,t\}$, then the path $(g,h_i)(g_1,h_1)(g_2,h_1)\cdots(t,h_1)(r,h_1)(t,h_2)\cdots(g'_2,h_2)(g'_1,h_2)(g',h_j)$ is the desired path. By symmetry, suppose that $g\in V(T)\backslash\{r,t\}$ and $g'\in\{r,t\}$. If $g'=r$, then the path $(g,h_i)(g_1,h_1)(g_2,h_1)\cdots(t,h_1)(r,h_j)$ is the desired path, while if $g'=t$, then the desired path is $(g,h_i)(g_1,h_1)(g_2,h_1)\cdots(t,h_1)(r,h_1)(t,h_j)$. Finally suppose that $g,g'\in\{r,t\}$. If $g\neq g'$, then the desired path is rather simple, that is the edge $(g,h_i)(g',h_j)$. For $g=g'=r$, if $i=1,j=2$ or $i=2,j\geq3$, then the path $(g,h_i)(t,h_2)(g',h_j)$ is the desired path; if $i=1,j\neq2$ or $i,j\geq3$, then the path $(g,h_i)(t,h_2)(a,h_2)(t,h_3)(g',h_j)$ is the desired path. For $g=g'=t$, if $i=1$, then the path $(g,h_i)(r,h_1)(g',h_j)$ is the desired path; if $i=2$ and $j\geq3$, then the path
$(g,h_i)(a,h_2)(g',h_j)$ is the desired path; if $3\leq i<j$, then the path $(g,h_i)(a,h_2)(t,h_2)(r,h_3)(g',h_j)$ is the desired path. Thus $tpc(T\circ H)\leq3$ and the proof is complete. \qed

\section{The strong product}

The {\it strong product} $G\boxtimes H$ of graphs $G$ and $H$ is the graph with vertex set $V(G)\times V(H)$, in which two vertices $(g,h),(g',h')$ are adjacent whenever $gg'\in E(G)$ and $h=h'$, or $g=g'$ and $hh'\in E(H)$, or $gg'\in E(G)$ and $hh'\in E(H)$. If an edge of $G\boxtimes H$ belongs to one of the first two types, then we call such an edge a {\it Cartesian edge} and an edge of the last type is called a {\it noncartesian edge}. (The name is due to the fact that if we consider only the first two types, we get the Cartesian product of graphs.) The strong product is commutative and is $2$-connected as long as both $G$ and $H$ are connected. Remind that $d_G(u,v)$ is the shortest distance between the two vertices $u$ and $v$ in graph $G$. And let $d_G(g)$ denote the degree of the vertex $g$ in $G$.

\begin{thm}\label{strong product} Let $G$ and $H$ be two nontrivial connected graphs. If $G\boxtimes H$ is not complete, then $tpc(G\boxtimes H)=3$.
\end{thm}

\pf Since $G\boxtimes H$ is not complete, $tpc(G\boxtimes H)\geq3$ and we only need to show that $tpc(G\boxtimes H)\leq3$. Like the method we use above, we pick a spanning tree $T$ of $G$ with root $t$ and a spanning tree $S$ of $H$ with root $s$. Clearly, the strong product $T\boxtimes S$ is a spanning subgraph of $G\boxtimes H$. Thus it suffices to show that $tpc(T\boxtimes S)\leq3$ by Proposition \ref{pro1}.

Define a total-coloring $c$ of $T\boxtimes S$ using the color set $A=\{1,2,3\}$ as follows. Let $W=\{w\in V(T)\backslash\{t\}: d_T(w)=1\}$ and $X=\{x\in V(S)\backslash\{s\}: d_S(x)=1\}$. We first give the vertices and edges of $T\boxtimes S$ a total-coloring using $A$ in such a way that for each $w\in W$, $x\in X$ and $v\in V(S)\backslash\{s\}$, the trail $(w,v)(p(w),v)\cdots(t,v)(t,p(v))\cdots(t,s)\cdots(w,s)\cdots(w,x)$ is total-proper connected and has the color sequence $1,3,2,1,3,2\cdots$ except for its two endpoints. Then set $c((u,s)(p(u),s^*))=c(p(u),s)$ where $u\in V(T)\backslash\{t\}$ and $s^*$ is a neighbor of $s$ in $S$. For any noncartesian edge, give it a color differing from the colors which its endpoints have already used. To complete our proof, the two claims below are necessary.

\textbf{Claim 1:} Let $tt_1t_2\cdots t_{i+j}$ and $ss_1s_2\cdots s_j$ be two paths in $T$ and in $S$ respectively, where $i\equiv0\ (\text{mod}\ 3)$ and $t_{i+j}\notin W$. Then for $x\in X$, the path $P=(t,x)(t,p(x))\cdots(t,s)(t_1,s)(t_2,s)\\\cdots(t_i,s)(t_{i+1},s_1)\cdots(t_{i+j},s_j)$ is total proper.

According to the total-coloring $c$ of $T\boxtimes S$, we only need to show that the path $P$ also has the color sequence $1,3,2,1,3,2\cdots$. If so, the vertices of $P$ should have the color sequence $1,2,3,1,2,3\cdots$. Consider the cycle $C=(t_{i+1},s_1)(t_i,s_1)\cdots(t,s_1)(t,s)(t_1,s)\cdots(t_i,s)(t_{i+1},s_1)$ and we have $|C|\equiv 0\ (\text{mod}\ 3)$. Then it follows that the cycle $C$ is total-proper connected and $c(t_{i+1},s_1)\equiv c(t_i,s)+1\ (\text{mod}\ 3)$. For the vertices $(t_{i+d},s_d)$ and $(t_{i+d+1},s_{d+1})$, where $1\leq d\leq j-1$, we set $P_1=(t_{i+d},s_d)(t_{i+d-1},s_d)\cdots(t,s_d)(t,s_{d-1})\cdots(t,s)$ and $P_2=(t_{i+d+1},s_{d+1})(t_{i+d},s_{d+1})\cdots(t,s_{d+1})(t,s_d)\cdots(t,s)$. Since $|P_2|=|P_1|+2$, we have $c(t_{i+d+1},s_{d+1})\equiv c(t_{i+d},s_d)-2\equiv c(t_{i+d},s_d)+1\ (\text{mod}\ 3)$. Thus, the vertices of the path $P$ do have the color sequence $1,2,3,1,2,3,\cdots$ and we complete the proof of Claim 1.

\textbf{Claim 2:} Let $t_it_{i+1}\cdots t_{i+k}$ and $s_js_{j+1}\cdots s_{j+k}$ be two paths in $T$ and in $S$ respectively, where $k\geq1,\ t_{i+k}\notin W,\ s_j\neq s,\ d_T(t_i,t)<d_T(t_{i+k},t)$ and $d_S(s_j,s)<d_S(s_{j+k},s)$. Then the path $P'=(t_i,s_j)(t_{i+1},s_{j+1})\cdots(t_{i+k},s_{j+k})$ is total-proper connected with the color sequence $1,3,2,1,3,2\cdots$.

Similar to the proof of Claim 1, for the vertices $(t_{i+d},s_d)$ and $(t_{i+d+1},s_{d+1})$, where $0\leq d\leq k-1$, we can deduce that $c(t_{i+d+1},s_{d+1})\equiv c(t_{i+d},s_d)-2\equiv c(t_{i+d},s_d)+1\ (\text{mod}\ 3)$. Then the vertices of $P'$ have the color sequence $1,2,3,1,2,3,\cdots$ and so the path $P'$ is total-proper connected with the color sequence $1,3,2,1,3,2\cdots$.

Next it remains to check that for any two vertices $(u,v),(u',v')$ in $T\boxtimes S$, where $u,u'\in V(T)$ and $v,v'\in V(S)$, there is a total proper path between them. Without loss of generality, assume that $d_S(v,s)\geq d_S(v',s)$. Let $tt_1t_2\cdots t_{a-1}u'$ and $ss_1s_2\cdots s_{b-1}v'$ be the paths from $t$ to $u'$ in $T$ and from $s$ to $v'$ in $S$, respectively. Then $d_T(t,u')=a$ and $d_S(s,v')=b$. If $a=b$, then the walk $(u,v)\cdots(t,v)\cdots(t,s)(t_1,s_1)(t_2,s_2)\cdots(u',v')$ is total-proper connected by Claim $1$, where $u\cdots t$ denotes the path from $u$ to $t$ in $T$ and $v\cdots s$ denotes the path from $v$ to $s$ in $S$. If $a<b$, then the walk $(u,v)(p(u),v)\cdots (t,v)(t,p(v))\cdots (t,s)(t_1,s_1)\\(t,s_1)(t_1,s_2)\cdots (t,s_{b-a})(t_1,s_{b-a+1})\cdots (t_{a-1},s_{b-1})(u',v')$ is total-proper connected by Claims $1$ and $2$. For $a>b$, assume that $u'$ lies on the path $tt_1t_2...t_{a-1}u'...w(=t_{d_T(w,t)})$ in $T$ from $t$ to $w$ where $w\in W$. If $v'=s$ or $u'=w$, then the situation is clear according to the total-coloring $c$ of $T\boxtimes S$. Otherwise we have $v'\neq s$ and $u'\neq w$. Then $b\geq1$ and $a\leq d_T(w,t)-1$. Let $c$ be a nonnegative integer. We divide the proof into three cases.

\textbf{Case 1.} $d_T(w,t)=3c+2$.

Then $a-b\leq 3c$ and so the walk $(u,v)(p(u),v)\cdots(t,v)(t,p(v))\cdots(t,s)(t_1,s)\cdots(t_{3d},s)\\(t_{3d+1},s_1)\cdots(t_{a-b+1},s_1)
(t_{a-b+2},s_2)\cdots(u',v')$ is total-proper connected by Claims $1$ and $2$, where $3(d-1)<a-b\leq 3d$ and $1\leq d\leq c$.

\textbf{Case 2.} $d_T(w,t)=3c+1$.

Then $a-b\leq 3c-1$. If $a-b\leq 3c-3$, we have the similar total-proper connected walk as Case 1. If $a-b=3c-1$, we must have $a=3c$ and $b=1$, and then the walk $(u,v)(p(u),v)\cdots(t,v)(t,p(v))\cdots(t,s)(t_1,s)\cdots(t_{a+1},s)(u',v')$ is total-proper connected. If $a-b=3c-2$, then we have $a=3c-1$ and $b=1$ or $a=3c$ and $b=2$. In the former case, we analogously have that the walk $(u,v)(p(u),v)\cdots(t,v)(t,p(v))\cdots(t,s)(t_1,s)\cdots(t_{a+1},s)(u',v')$ is total-proper connected; while in the latter case, since $d_T(w,t)+b=3c+3\equiv 0\ (\text{mod}\ 3)$, the walk
$(u,v)(p(u),v)\cdots(t,v)(t,p(v))\cdots(t,s)(t_1,s)\cdots(w,s)(w,s_1)(w,v')(u',v')$ is total-proper connected.

\textbf{Case 3.} $d_T(w,t)=3c$.

Then $a-b\leq 3c-2$. If $a-b\leq 3c-3$, we have the similar total-proper connected walk as Case 1. Thus $a=3c-1$ and $b=1$ is the only case we need to discuss. However this is like the discussion in Case 2.

After removing all the cycles of the total-proper connected walks appearing above, we get the corresponding desired paths. Thus, $tpc(T\boxtimes S)\leq3$ and we complete the proof.\qed


\begin{thebibliography}{1}

\bibitem{ACKP}
 B. S. Anand, M. Changat, S. Klav$\check{z}$ar, I. Peterin, Convex sets in lexicographic products of graphs, {\it Graphs Combin.} {\bf 28(1)} (2012) 77-84.

\bibitem{ALL}
 E. Andrews, E. Laforge, C. Lumduanhom, P. Zhang, On proper-path colorings
 in graphs, {\it J. Combin. Math. Combin. Comput} {\bf 97} (2016) 189-207.

\bibitem{B}
 J.A. Bondy, U.S.R. Murty, {\it Graph Theory}, GTM 244, Springer, 2008.

\bibitem{BFG}
 V. Borozan, S. Fujita, A. Gerek, C. Magnant, Y. Manoussakis, L. Montero,
 Z. Tuza, Proper connection of graphs, {\it Discrete Math.} {\bf 312(17)} (2012) 2550-2560.

\bibitem{BS}
 B. Bre$\check{s}$ar, S. $\check{S}$pacapan, On the connectivity of the direct product of graphs, {\it Australas. J. Comb.} {\bf 41} (2008) 45-56.

\bibitem{CH}
 G. Chartrand, F. Harary, Planar permutation graphs, {\it Ann. Inst. H. Poincar$\acute{e}$ (Sect. B)} {\bf 3(4)}(1967), 433-438.

\bibitem{GLQ}
 R. Gu, X. Li, Z. Qin, Proper connection number of random graphs, {\it Theoret. Comput. Sci.} {\bf 609(2)} (2016), 336--343.

\bibitem{GMP}
 T. Gologranc, G. Meki$\check{s}$, I. Peterin, Rainbow connection and graph products, {\it Graphs Combin.} {\bf 30(3)} (2014) 591-607.

\bibitem{GV}
 R. Guji, E. Vumar, A note on the connectivity of Kronecker products of graphs, {\it Appl. Math. Lett.} {\bf 22(9)} (2009) 1360-1363.

\bibitem{JLZ}
 H. Jiang, X. Li, Y. Zhang, Upper bounds for the total rainbow connection of graphs, {\it J. Comb. Optim.} {\bf 32(1)} (2016) 260-266.

\bibitem{JLZ1}
 H. Jiang, X. Li, Y. Zhang, Total proper connection of graphs, arXiv:1512.00726v1.

\bibitem{JLZZ}
 H. Jiang, X. Li, Y. Zhang, Y. Zhao, On (strong) proper vertex-connection of graphs, {\it Bull. Malays. Math. Sci. Soc.} {\bf 93(1)} (2015) 1-11.

\bibitem{KS}
 S. Klav$\check{z}$ar, S. $\check{S}$pacapan, On the edge-connectivity of Cartesian product graphs, {\it Asian-Eur. J. Math.} {\bf 1} (2008) 93-98.

\bibitem{LLZ}
 E. Laforge, C. Lumduanhom, P. Zhang, Characterizations of graphs
 having large proper connection numbers, {\it Discuss. Math. Graph Theory} {\bf
36(2)} (2016) 439-453.

\bibitem{LiSS}
 X. Li, Y. Shi, Y. Sun, Rainbow connections of graphs: A survey, Graphs \& Combin. {\bf 29(1)} (2013), 1-38.

\bibitem{LiS}
 X. Li, Y. Sun, Rainbow Connections of Graphs, Springer Briefs in Math., Springer, New York, 2012.

\bibitem{LMS}
 H. Liu, $\hat{A}$. Mestre, T. Sousa, Total rainbow $k$-connection in graphs, {\it Discrete Appl. Math.} {\bf 174} (2014), 92-101.

\bibitem{MYWY}
 Y. Mao, F. Yanling, Z. Wang, C. Ye, Rainbow vertex-connection and graph products, {\it Int. J. Comput. Math.} {\bf 93(7)} (2016) 1078-1092.

\bibitem{MYWY1}
 Y. Mao, F. Yanling, Z. Wang, C. Ye, Proper connection number and graph products, {\it Bull. Malays. Math. Sci. Soc.} DOI 10.1007/s40840-016-0442-z, in press.

\bibitem{NS}
 R. J. Nowakowski, K. Seyffarth, Small cycle double covers of products. I. Lexicographic product with paths and cycles, {\it J. Graph Theory} {\bf 57(2)} (2008) 99-123.

\bibitem{P}
 I. Peterin, Intervals and convex sets in strong product of graphs, {\it Graphs Combin.} {\bf 29(3)} (2013) 705-714.

\bibitem{S}
 Y. Sun, On rainbow total-coloring of a graph, {\it Discrete Appl. Math.} {\bf 194}(2015), 171-177.

\bibitem{S1}
 S. $\check{S}$pacapan, Connectivity of strong products of graphs, {\it Graphs Combin.} {\bf 26(3)} (2010) 457-467.

\bibitem{Z}
 X. Zhu, Game coloring the Cartesian product of graphs, {\it J. Graph Theory} {\bf 59(4)} (2008) 261-278.
\end{thebibliography}
\end{document}